\documentclass[11pt,a4paper,notitlepage,final,onecolumn,oneside]{article}
\usepackage{exscale}
\usepackage{latexsym}
\usepackage{calc}
\usepackage{amssymb}
\usepackage{url}

\newcommand{\proof}{{\bfseries Proof:\quad}}
\newcommand{\p}{\ensuremath{\hat{\mathrm{p}}}}

\begin{document}

\newtheorem{lemma}{Lemma}
\newtheorem{proposition}[lemma]{Proposition}
\newtheorem{theorem}[lemma]{Theorem}
\newtheorem{definition}[lemma]{Definition}
\newtheorem{hypothesis}[lemma]{Hypothesis}
\newtheorem{conjecture}[lemma]{Conjecture}
\newtheorem{remark}[lemma]{Remark}
\newtheorem{example}[lemma]{Example}
\newtheorem{property}[lemma]{Property}
\newtheorem{corollary}[lemma]{Corollary}
\newtheorem{algorithm}[lemma]{Algorithm}


\title{The $k$-orbit theory and Fein-Kantor-Schacher Theorem}
\author{Aleksandr Golubchik\\
        \small Bramscher Str. 57, 49088 Osnabrueck, Germany,\\
        \small e-mail: agolubch@uos.de}

\date{02.02.2005}
\maketitle


\begin{abstract}
By the investigation of $k$-orbits symmetry properties it is obtained a
simple proof of the B. Fein, W. M. Kantor and M. Schacher Theorem: any
transitive permutation group contains a non-trivial fixed-point-free
prime-power element. \bigskip


\end{abstract}


\section{Introduction}

Using the Classification of finite simple groups, B. Fein, W. M. Kantor
and M. Schacher \cite{FKS} give a proof of the

\begin{theorem}\label{FKS}
Any transitive permutation group contains a non-trivial fixed-point-free
prime-power element.

\end{theorem}

In the pre-doctoral and doctoral students course ``Permutation groups,
structures, and polynomials'' in Prague, January - February 2004, P.
Cameron suggested a problem: to find a simple proof of the theorem
\ref{FKS} that does not use the Classification of finite simple groups.
This problem P. Cameron estimated as ``very difficult''.

Below we give such a proof, based on the $k$-orbit theory developed by
author that studies symmetry properties of $k$-orbits. The first attempt
to describe this theory was launched in \cite{Golubchik-2}.


\section{$k$-Orbit theory}

Let $G(V)$ be a permutation group of a degree $n$ and $V^{(k)}$ ($k\leq
n$) be a non-diagonal part of $V^k$ (i.e. all $k$ values of coordinates of
a $k$-tuple $\alpha_k\in V^{(k)}$ are different), then the action of $G$
on $V^{(k)}$ forms a partition $GV^{(k)}$ of $V^{(k)}$ on $G$-invariant
classes. This partition is called a system of $k$-orbits of $G$ and we
denote it as $Orb_k(G)$.

The $k$-orbit theory studies symmetry properties of $k$-orbits that cannot
be obtained from permutation group theory, studying a permutation group
$G$ as a permutation algebra and therefore being not able to see the
inside structure of a $n$-orbit of $G$ and relations between $k$-orbits of
subgroups of $G$.

The $k$-orbit theory gives a new view on abstract finite groups and finite
permutation groups and shows a way to simple solutions of some problems
which either are not solved or have a complicated solution. These problems
are: the polynomial solution of graphs isomorphism problem; the problem of
a full invariant of finite groups, admitting a primitive (isomorphic)
permutation representation; Polycirculant conjecture: the automorphism
group of a transitive graph contains a regular element (s.
\cite{Cameron}); a simple proof of the Feit-Thompson theorem: the
solvability of finite groups of odd order (that lies in the foundation of
the Classification of finite simple groups \cite{Gorenstein}) and the
simplest of those task that we consider in this paper.

The specificity of the $n$-orbit representation of a group $G$ is a
possibility to do a group visible. In order to go to this visibility one
makes a partition of a matrix of a $n$-orbit on cells of $k$-orbits of
subgroups of an investigated permutation group $G$ and studies symmetry
properties of cells and the whole partition. Below we consider some facts
from the $k$-orbit theory which are necessary for the proof of theorem
\ref{FKS}.

The $k$-orbit theory considers different actions of permutations on
$k$-sets and symmetries that follows from those actions. Let
$\alpha_n=\langle v_1\ldots v_n\rangle$ be a $n$-tuple, then the left
action of a permutation $g$ on $\alpha_n$ is defined as $g\alpha_n=\langle
gv_1\ldots gv_n\rangle$ and the right action of $g$ on $\alpha_n$ is given
as $\alpha_ng=\langle v_{g1}\ldots v_{gn}\rangle$ (here $V=[1,n]$). From
that we obtain the left and right action of a permutation on $k$-tuples
and $k$-sets.

The left action of a permutation is (on definition) an isomorphism
(\,$(123)\{123,132\}=\{231,213\}$\,), the right action is not an
isomorphism (\,$\{123,132\}(123)=\{231,321\}$\,). If the right action is an
isomorphism, then it is an indicium for existing of a subgroup with
non-trivial normalizer, whose $n$-orbit contains (right-)isomorphic
$k$-subsets connected with considered right action, for example a subgroup
$$
\begin{array}{||cc||}
\hline
12 & 34 \\
21 & 43 \\
\hline
\end{array}\
\mbox{ in a group }\
\begin{array}{||cc||}
\hline
12 & 34 \\
21 & 43 \\
\hline
34 & 12 \\
43 & 21 \\
\hline
\end{array}\,,
$$
where subgroup and group are represented through their $n$-orbits.

Let $X_n$ be a $n$-orbit of $G$. One of $n$-tuples from $X_n$ we choose as
an \emph{initial} $n$-tuple and consider all permutations from $G$
relatively to this $n$-tuple. If $\alpha_n$ is an initial $n$-tuple, then
$X_n=G\alpha_n$. The specificity of an initial $n$-tuple is the equality
of number and order values of its coordinates (that we take from the same
ordered set $V$). We indicate the $n$-space $V^n$ through the initial
$n$-tuple. So the ordered subsets of the initial $n$-tuple point out
corresponding subspaces. A $k$-orbit $X_k$ of $G$ is a projection of $X_n$
on some subspace $I_k\in X_k$. We write this projection as
$X_k=\p^{(I_k)}X_n$. So we consider subspaces as ordered sets of
coordinates and $k$-orbits as non-ordered sets of $k$-tuples. The
different orderings of $k$-tuples (lines) in a matrix of a $k$-orbit shows
different symmetry properties of $k$-orbit, related with corresponding
properties of the investigated permutation group. We can consider of
course also non-ordered by coordinates $k$-orbits or subspaces, if it is
suitable. In such case a $k$-orbit is represented by matrix accurate to
permutation of columns.

We say that two $k$-sets $X_k$ and $X_k'$ are $G$-isomorphic if they are
isomorphic and connected with a permutation $g\in G$:  $X_k'=gX_k\equiv
\{g\alpha_k:\,\alpha_k\in X_k\}$.  If we study a group $G$ and $X_k,X_k'$
are isomorphic, but not $G$-isomorphic, then we say that they are
$S_n$-isomorphic (where $S_n$ is the symmetric group).

We defined the left action $GY_n$ of a group $G$ on a $n$-subset
$Y_n$ as $GY_n\equiv \{gY_n:\, g\in G\}$. The left action $AX_n$ of a
subgroup $A<G$ on a $n$-orbit $X_n$ of $G$ is given as $AX_n\equiv
\{A\alpha_n:  \alpha_n\in X_n\}$. Correspondingly the right action
$Y_nG\equiv \{Y_ng:\, g\in G\}$ and $X_nA\equiv \{\alpha_nA:\alpha_n\in
X_n\}$.  The left action of $G$ on a $k$-set $Y_k$ and $A$ on a $k$-orbit
$X_k$ of $G$ for $k<n$ is the same as for $k=n$, because
$\p^{(I_k)}g\alpha_n=g\p^{(I_k)}\alpha_n$.  By the right action the
equality $\p^{(I_k)}(\alpha_ng)= (\p^{(I_k)}\alpha_n)g$ has sense only if
$\alpha_n$ is the initial $n$-tuple. Hence the right action on $k$-sets
has not a direct reduction from the right  action on $n$ sets. If $X_n$ is
a $n$-orbit of $G$, $A<G$ and $Y_n$ is a $n$-orbit of $A$, then
$L_n=GY_n=X_nA$ and $R_n=Y_nG=AX_n$ are partitions of $X_n$ on $n$-orbits
of left and right cosets of $A$ in $G$, at that $n$-orbits of left cosets
are $n$-orbits of conjugate to $A$ subgroups of $G$, because $n$-orbits of
left cosets are $G$-isomorphic, $n$-orbits of right cosets are $n$-orbits
of $A$ , because $n$-orbits of right cosets differ in order of
coordinates. Projections of $X_n$, $L_n$ and $R_n$ on a subspace $I_k$
gives a covering $L_k=GY_k$, where $Y_k=\p^{(I_k)}Y_n$, and a partition
$R_k=AX_k$ of the $k$-orbit $X_k$ on $k$-orbits of left and right cosets
of $A$ in $G$. Below under an action of permutations we understand the left
action.

Let $k$ be a divisor of $|G|$. We call $k$ as an \emph{automorphic number}
if there exists a $k$-element \emph{suborbit} of $G$ (a $k$-element orbit
of a subgroup of $G$).

Let $Co(\alpha_k)$ be a set of coordinates of a $k$-tuple $\alpha_k$, then
for a $k$-set $X_k$ we define $Co(X_k)\equiv\{Co(\alpha_k):\,\alpha_k\in
X_k\}$.

If $Co(\alpha_k)$ is a suborbit of $G$, then we say that $\alpha_k$ is an
\emph{automorphic $k$-tuple} and $Co(\alpha_k)$ is \emph{an automorphic
subset} of $V$.  Let $X_k$ be a $k$-set we say that $X_k$ is an
\emph{automorphic $k$-set}, if it is a $k$-orbit. Let $X_k$ be a $k$-orbit
and $\alpha_k\in X_k$ be an automorphic $k$-tuple, then a $k$-orbit $X_k$
we call as \emph{right-automorphic $k$-orbit} or \emph{$k$-rorbit}.

Let $U$ be a set and $Q$ be a set of subsets of $U$, then $\cup Q\equiv
\cup_{(W\in Q)}W$. Under $\sqcup Q$ we understand a set of unions of
intersected on $U$ classes of $Q$. So $\sqcup Q$ is a partition of $\cup
Q$ (possibly trivial). A union and intersection of partitions $P$ and $R$
of $U$ we write as $P\sqcup R$ and $P\sqcap R$. Let $P'$ be a subpartition
of $P$, then we write $P'\sqsubset P$.

Let $X_k$ be a $k$-rorbit, $\alpha_k\in X_k$ and $Y_k$ be a maximal subset
of $X_k$ so that $Co(Y_k)=Co(\alpha_k)$, then we call $Y_k$ as a
\emph{$k$-block} of $X_k$. If $Y_k$ is a $k$-block, then $Aut(Y_k)$ is a
transitive group of degree $k$.

Let $X_k$ be a $k$-rorbit, $\cup Co(X_k)=V$ and $\sqcup Co(X_k)$ be a
non-trivial partition of $V$, then we say that $X_k$ is \emph{incoherent}.
If $\sqcup Co(X_k)$ is a non-trivial covering of $V$, then we say that
$X_k$ is \emph{coherent}. Any coherent $k$-orbit consists of intersected
on $V$ $k$-blocks. An incoherent $k$-orbit contains $k$-suborbits that
form a partition of the $k$-orbit, where classes of this partition are
non-intersected on $V$ coherent $k$-suborbits or $k$-blocks. A coherent
$k$-orbit $X_k$ is defined on a set $V=\cup Co(X_k)$. Its coherent or
incoherent $k$-suborbit $Y_k$ is defined on a set $U=\cup Co(Y_k)\subset
V$. In order to show this we say that $Y_k$ is $U$-coherent
($U$-incoherent). If a ($V$-) coherent $k$-orbit contains no $U$-coherent
and no $U$-incoherent $k$-suborbit for any $U\subset V$, then we call it
as \emph{elementary coherent}.

Let $Y_k$ be a $k$-suborbit of $G$. The maximal transitive on $Y_k$
subgroup of $G$ we call a \emph{stabilizer} of $Y_k$ in $G$ and write it as
$Stab_{\, G}(Y_k)$ or simply $Stab(Y_k)$.

Let $X_k$ be a $k$-orbit of $G$ and $Y_k\subset X_k$ be a $k$-block, then
evidently $L_k=GY_k$ is a partition of $X_k$ and $Y_k$ is a $k$-orbit of
$Stab(Co(Y_k))$.

Let a group $G$ be imprimitive, then $G$ contains a non-trivial ($1<k<n$)
incoherent $k$-orbit. If $G$ is primitive, then any (non-trivial)
$k$-orbit of $G$ is coherent. From here follows that it is convenient to
consider primitive Abelian groups as trivial imprimitive, because, as
distinct from non-Abelian primitive groups, such group $G$ contains no
non-trivial suborbit $U\subset V$ that forms a covering $GU$ of $V$. So
further under primitive group we understand non-Abelian primitive group.
Let a group $G$ be imprimitive and $Q$ be a $G$-invariant partition of $V$
($GQ=Q$), then we say that classes of $Q$ are \emph{imprimitivity blocks}.


\subsection{$k$-orbits and normal subgroups}

$k$-Orbit theory shows a new approach for consideration of normal
subgroups that based on pure combinatorial symmetries of $k$-orbits. We
consider here some statements that show general relations between
$k$-orbits of normal subgroups and statements that we need for the proof
of theorem \ref{FKS}.

\begin{proposition}\label{Stab(Co(a_k))}
Let $\alpha_k$ be an automorphic $k$-tuple, $G=Stab(Co(\alpha_k))$,
$H=Stab(\alpha_k)$ and $X_k=G\alpha_k$, then $H$ is a normal subgroup
of $G$ and a factor group $G/H$ is isomorphic to $Aut(X_k)$.

\end{proposition}
\proof
If $H$ is not trivial, then $G$ is an intransitive group. Let $X_n$ be a
$n$-orbit of $G$, then there exists only one partition $P_n$ of $X_n$ on
$|G|/|H|$ classes with projection $\p^{(\alpha_k)}P_n=X_k$. Hence
$X_nH=HX_n=P_n$. The group $Aut(X_k)$ is isomorphic to the factor group
$G/H$, because $|Aut(X_k)|=|G/H|$ and it acts on classes of $H$ in $G$
transitive. $\Box$\bigskip

\begin{proposition}\label{L_k=R_k}
Let $X_k\in Orb_k(G)$, $Y_k\subset X_k$ be a $k$-suborbit of $G$,
$H=Stab(Y_k)$, $L_k=GY_k$ and $R_k=HX_k$. Let $L_k=R_k$, then $H$ is a
normal subgroup of $G$.

\end{proposition}
\proof
There exists only one partition $P_n$ of $X_n$ on $|G|/|H|$ classes so
that $\p^{(I_k)}P_n=L_k$. Hence $X_nH=HX_n=P_n$. $\Box$\bigskip

\begin{proposition}\label{G.pr-H.tr}
Let $H$ be a normal subgroup of a primitive group $G$, then $H$ is
transitive.

\end{proposition}
\proof
Let $H$ be intransitive normal subgroup of a transitive group $G$, then
projections of $Y_n=HI_n$ on orbits of $H$ are $G$-isomorphic and hence a
normalizer $N_G(H)$ is imprimitive. $\Box$\bigskip

\begin{corollary}\label{simple}
Let $G$ be a primitive group with no transitive subgroup, then $G$ is a
simple group.

\end{corollary}

\begin{proposition}\label{A<H<|G}
Let $A<H\lhd G$ be abstract groups and $N_H(A)=A$, then $N_G(A)\neq A$ and
$|G|/|H|=|N_G(A)|/|A|$.

\end{proposition}
\proof
$A^G=A^H$ and $|A^H|=|H|/|A|$. so
$|N_G(A)|/|A|=(|G|/|A|)/(|H|/|A|)=|G|/|H|$. $\Box$\bigskip

\begin{proposition}\label{A(V)<H(V)<|G(V)}
Let $A(V)<H(V)\lhd G(V)$ be permutation groups and $N_H(A)=A$, then $A$
contains $G$-isomorphic $k$-orbits.

\end{proposition}
\proof
First, the statement follows from proposition \ref{A<H<|G}.
Second, $H$ has $G$-isomorphic $k$-orbits $Y_k'=gY_k$ (as a normal
subgroup). So $Y_k'$ has the same structure as $Y_k$, i.e. if $Y_k$
contains a $k$-orbit $Z_k$ of $A$, then $Y_k'$ contains also a $k$-orbit
$Z_k'$ of $A$ and $Z_k'=gZ_k$ for some $g\in G$. $\Box$\bigskip

\begin{lemma}\label{G<A_n}
Let $G<A_n$, then $N_{S_n}(A)\neq A$.

\end{lemma}

\begin{corollary}\label{no.tr.subgr.}
Let $G$ be a primitive group with no transitive subgroup, then
$N_{S_n}(G)\neq G$.

\end{corollary}


\subsection{Some relations between $|X_k|$ and $|Aut(X_k)|$}

The consideration of a $k$-orbit $X_k$ is easier if $|X_k|=|Aut(X_k)|$.
Here we consider some conditions for this equality.

\begin{proposition}\label{|Aut(X_k)|=|X_k|}
Let $X_k$ be a $k$-orbit and for every subgroup $A<Aut(X_k)$ a partition
$R_k=AX_k$ consists of classes of the same power, then $|Aut(X_k)|=|X_k|$.

\end{proposition}
\proof
In such case no permutation from $Aut(X_k)$ fixes a $k$-tuple from $X_k$.
$\Box$\bigskip

\begin{lemma}\label{|gr(A,B)|}
Let $X_k$ be a $k$-orbit, $Y_k,Z_k\subset X_k$ be $k$-suborbits that have
non-trivial intersection, $\alpha_k \in Y_k\cap Z_k$, $A,B<Aut(X_k)$,
$Y_k=A\alpha_k$, $Z_k=B\alpha_k$, $|Y_k|=|A|$ and $|Z_k|=|B|$. Let
$T_k\subset X_k$ be a $k$-orbit of a group $G=gr(A,B)$, then $|T_k|=|G|$.

\end{lemma}
\proof
Let $Y_n\in Orb_n(A)$ and $Z_n\in Orb_n(B)$ , then $|Y_k|=|Y_n|$ and
$|Z_k|=|Z_n|$. Let $T_n\in Orb_n(G)$, then $|T_n|/|Y_n|=|T_n|/|Y_k|$ and
$|T_n|/|Z_n|=|T_n|/|Z_k|$. It follows that $L_k'=GY_k$ and $L_k''=GZ_k$
are partitions of $T_k=L_k'\sqcup L_k''$. $\Box$\bigskip

\begin{proposition}\label{cap,cup}
Let $X_k\in Orb_k(G)$, $Y_k,Z_k\subset X_k$ be $k$-suborbits of $G$ and
$L_k'=GY_k,L_k''=GZ_k$ be partitions of $X_k$, then $L_k^1=L_k'\sqcap
L_k''$ and $L_k^2=L_k'\sqcup L_k''$ are also partitions of $X_k$ on
isomorphic $k$-suborbits of $G$. If $Y_k\cap Z_k=T_k$ is not trivial and
$U_k\in L_k^2$ contains $Y_k$ and $Z_k$, then $L_k^1=GT_k$, $L_k^2=GU_k$,
$Stab(T_k)=Stab(Y_k)\cap Stab(Z_k)$ and
$Stab(U_k)=gr(Stab(Y_k),Stab(Z_k))$.

\end{proposition}
\proof
Let $k=n$, then equalities follows from corresponding properties of sets of
left cosets of subgroups of $G$. For $k<n$ equalities are projections of
corresponding equalities for $k=n$. $\Box$\bigskip

\begin{lemma}\label{|H|=|X_k|}
Let for every $k$-suborbit $Y_k$ of a $k$-orbit $X_k$ a set
$L_k=Aut(X_k)Y_k$ be a partition of $X_k$, then $Aut(X_k)$ contains a
transitive on $X_k$ normal subgroup $H$ of order $|X_k|$.

\end{lemma}
\proof
Let $G=Aut(X_k)$, $M$ be a set of partitions of $X_k$ on isomorphic
$k$-suborbits of $G$ and $L_k',L_k''\in M$, then $L_k^1=L_k'\sqcap L_k''$
and $L_k^2=L_k'\sqcup L_k''$ are partitions from $M$. Let $L_k''$ be a
partition of $X_k$ on $k$-blocks, then there exists a partition $L_k'$ so
that $L_k^2=L_k'\sqcup L_k''$ is not trivial, i.e.
$|L_k^2|<min(|L_k'|,|L_k''|)$. We can assume now $L_k''=L_k^2$ and then
obtain new $L_k^2=L_k'\sqcup L_k''$ for corresponding new partition
$L_k'$. By repeating we obtain that there exist partitions $L_k',L_k''\in
M$ so that $L_k^2=L_k'\sqcup L_k''=\{X_k\}$.

We can realize a such union process for that a class $Y_k$ of a partition
$L_k$ is a $k$-orbit of a subgroup $A<G$ of order $|A|=|Y_k|$, because we
can begin from a subgroup that is isomorphic to the automorphism group of
a $k$-block, and with a subgroup of a prime order that permutes
$k$-blocks. Thus, according to lemma \ref{|gr(A,B)|}, $G$ contains a
subgroup $H$ of order $|X_k|$ that acts transitive on $X_k$. This subgroup
is normal in $G$, because there exists only one partition of a $n$-orbit
$X_n$ of $G$ on automorphic classes with the same projection $X_k$ on a
subspace from $X_k$. $\Box$\bigskip

\begin{proposition}\label{}
Let $X_k$ be incoherent, then $|Aut(X_k)|\neq |X_k|$.

\end{proposition}
\proof
Let $m$ be a number of $k$-blocks of $X_k$ and $Y_k$ be a $k$-block. If
$m>2$, then $X_k$ contains a $k$-suborbit of a power $(m-1)|Y_k|$ that
does not divide $|X_k|$. If $m=2$, then
$|Stab(Y_k)|=|Y_k|^2$.~$\Box$\bigskip

\begin{proposition}\label{N_G(A)=A}
Let in lemma \ref{|H|=|X_k|} the group $G$ contains a subgroup $A$ with
trivial normalizer, then $H=G$, i.e. $|Aut(X_k)|=|X_k|$.

\end{proposition}
\proof
Let $R_k=AX_k$ and $Y_k\in R_k$, then $Stab(Y_k)=A$. Let $H\neq G$, then
$L_k=GY_k$ is a covering of $X_k$.  Contradiction. Hence $H=G$.
$\Box$\bigskip

\begin{theorem}\label{X_k-coh.}
Let in lemma \ref{|H|=|X_k|} $X_k$ be coherent, then $|Aut(X_k)|=|X_k|$.

\end{theorem}
\proof
Let $A=Stab(v)$, $v\in V$ and $R_k=AX_k$. The $k$-orbit $X_k$ defines a
group $G=Aut(X_k)$ that is generated by its stabilizer $A$ and other
subgroups which belong to no stabilizer. Hence those generators have to be
defined by $k$-suborbits from $X_k$. It follows that $R_k$ contains a class
$Y_k$ with a power $|Y_k|=|A|$. So if $H\neq G$, then $L_k=GY_k$ is a
covering of $X_k$.  Contradiction. $\Box$\bigskip

\begin{lemma}\label{Aut(el.coh)}
Let $X_k$ be an elementary coherent $k$-orbit, then for every $k$-suborbit
$Y_k\subset X_k$ a set $L_k=Aut(X_k)Y_k$ is a partition of $X_k$.

\end{lemma}
\proof
Let $Y_k\subset X_k$ be a $k$-block, then $L_k=GY_K$ is a partition of
$X_k$ on definition. Let $Z_k\subset X_k$ be a $k$-suborbit
of $Aut(X_k)$ that has non-trivial intersection with $Y_k$. Let
$A=gr(Stab(Y_k),Stab(Z_k))$, then $k$-orbit of $A$ in $X_k$ is $X_k$,
because the elementary coherent $k$-orbit $X_k$ is the unique
super-$k$-suborbit for any $k$-block. Hence $Aut(X_k)Z_k$ is a partition
of $X_k$. $\Box$\bigskip

\begin{theorem}\label{el.coh.-H}
Let $X_k$ be an elementary coherent $k$-orbit, then $|Aut(X_k)|=|X_k|$.

\end{theorem}


\section{A simple Proof of the Theorem \ref{FKS}}

Let $G(V)$ be an imprimitive group and $Q$ be a partition of $V$ on
imprimitivity blocks, then $G$ has a homomorphic (possibly isomorphic)
transitive representation $G'(Q)$ on classes of $Q$.

\begin{lemma}\label{G'(Q)}
Let $G'(Q)$ contains a fixed-point-free element of a prime-power order,
then $G(V)$ contains a fixed-point-free element of a prime-power order too.

\end{lemma}
\proof
Let $g'\in G'$ be a fixed-point-free element of a prime-power order $p^l$
that is a reduction of an element $g\in G$, then $g$ is evidently a
fixed-point-free element and $|g|=p^md$, where a multiplier $d$ is coprime
to $p$ and $m\geq l$. So $g^d$ is a fixed-point-free element of a
prime-power order of the group $G$. $\Box$\bigskip

Hence it is sufficient to consider the theorem \ref{FKS} for primitive
groups.

Let $G(V)$ be a primitive group that contains a transitive subgroup $A$,
then we can reduce the task on subgroup $A$. So it is necessary to
consider the theorem \ref{FKS} for a primitive group that contains no
transitive subgroup. Let $G(V)$ be a such group of a degree $n$ and $k$ be
a maximal automorphic divisor of $n$ .

We know that a normalizer $N=N_{S_n}(G)\neq G$, i.e. $G$ contains
$N$-isomorphic $k$-orbits. Hence there must exist a partition
$Q$ of $V$ on $k$-element $N$-isomorphic suborbits of $G$, so that
projections of a $n$-orbit $X_n$ of $G$ on these suborbits are
$N$-isomorphic $k$-orbits.

\begin{lemma}\label{}
Let $X_k$ be a such $k$-orbit, then it is elementary coherent.

\end{lemma}
\proof
Let $X_k$ be not elementary coherent, then there exists a suborbit
$U\subset V$ of $G$ that consists of not intersected $k$-tuples from
$X_k$ and $|U|$ does not divide $|V|$. It gives a contradiction in action
of $N$ on $U$ and on the partition $Q$. $\Box$\bigskip

Since $X_k$ is elementary coherent, $|X_k|=|G|$. From here follows that a
stabilizer $A$ of a class from $Q$ is isomorphic to its projections on
classes of $Q$. If a such projection of $A$ has fixed-point-free
prime-power element, then $A$ also has fixed-point-free prime-power
element. Thus we have again a reduction of a task that solves theorem
\ref{FKS}.


\section*{Conclusion}

The main inference from the investigation of $k$-orbits symmetry
properties is that a finite group $F$ is not closed by its algebraic
properties, because the group algebra $F$ is equivalent to the action of
$F$ on $F$, but this algebra generates also the action of $F$ on $F^k$.
Properties of that action not always can be interpreted with group algebra
or with traditional permutation group theory. Namely such properties are
the subject of investigation of the $k$-orbit theory. Some investigations
related with an application of $k$-orbit representation to problems
announced above are described in \cite{Golubchik-1} - \cite{Golubchik-3}

The $k$-orbit theory gives a new view on finite many-dimensional
symmetries and can bring new ideas and applications to other sciences
studying symmetrical objects.


\section*{Acknowledgements}
I would like to express many thanks to Prof.~P.~Cameron for his excellent
site on the internet, for two years of e-mail contacts and for presenting
me some examples of permutation groups interesting for analysis by methods
of $k$-orbit theory. I also would like to express many thanks to Prof.
Bernd Fischer for his patience to hear me out.


\end{document}